\documentclass[journal,11pt,draftcls,onecolumn,peerreview]{IEEEtran}

\usepackage{graphicx}
\usepackage{url}
\usepackage{amsmath,amssymb,amsthm,bm}
\usepackage{subfig}
\usepackage{cite}
\usepackage[latin1]{inputenc}

\newtheorem{thm}{Lemma}[subsection]

\newcommand{\abs}[1]{\lvert#1\rvert}
\newcommand{\edmc}{\ensuremath{\cal{E}}}

\begin{document}

\title{Robust Simultaneous Localization of Nodes and Targets in Sensor Networks Using Range-Only Measurements}

\author{P\i nar~O\u{g}uz-Ekim,~\IEEEmembership{Student Member,~IEEE,}
        Jo\~{a}o~Gomes,~\IEEEmembership{Member,~IEEE,}
        Jo\~{a}o~Xavier,~\IEEEmembership{Member,~IEEE,}
        and~Paulo~Oliveira,~\IEEEmembership{Member,~IEEE,}% <-this % stops a space
\thanks{The authors are with the Institute for Systems and Robotics -- Instituto Superior T\'ecnico (ISR/IST), Lisbon, Portugal. e-mail: poguz@isr.ist.utl.pt.}% <-this % stops a space
\thanks{This research was partially supported by Funda\c{c}\~{a}o para a Ci\^{e}ncia e a Tecnologia (FCT) through ISR/IST plurianual funding with PIDDAC program funds, project PTDC/EEA-CRO/104243/2008, and grant SFRH/BD/44771/2008.}}

\ifCLASSOPTIONpeerreview
\else
\markboth{IEEE TRANSACTIONS ON SIGNAL PROCESSING}
{O\u{g}uz-Ekim \MakeLowercase{\textit{et al.}}: Robust simultaneous localization of nodes and targets in sensor networks using range-only measurements}
\fi

\maketitle

\begin{abstract}

Simultaneous localization and tracking (SLAT) in sensor networks aims to determine the positions of sensor nodes and a moving target in a network, given incomplete and inaccurate range measurements between the target and each of the sensors. One of the established methods for achieving this is to iteratively maximize a likelihood function (ML), which requires initialization with an approximate solution to avoid convergence towards local extrema. This paper develops methods for handling both Gaussian and Laplacian noise, the latter modeling the presence of outliers in some practical ranging systems that adversely affect the performance of localization algorithms designed for Gaussian noise. A modified Euclidean Distance Matrix (EDM) completion problem is solved for a block of target range measurements to approximately set up initial sensor/target positions, and the likelihood function is then iteratively refined through Majorization-Minimization (MM). To avoid the computational burden of repeatedly solving increasingly large EDM problems in time-recursive operation an incremental scheme is exploited whereby a new target/node position is estimated from previously available node/target locations to set up the iterative ML initial point for the full spatial configuration. The above methods are first derived under Gaussian noise assumptions, and modifications for Laplacian noise are then considered. Analytically, the main challenges to be overcome in the Laplacian case stem from the non-differentiability of $\ell_1$ norms that arise in the various cost functions. Simulation results confirm that the proposed algorithms significantly outperform existing methods for SLAT in the presence of outliers, while offering comparable performance for Gaussian noise.
\end{abstract}

\ifCLASSOPTIONpeerreview
% \begin{center} \bfseries EDICS Category: SEN-LOCL \end{center}

% \noindent\textbf{Contact Information:} \\
% Institute for Systems and Robotics --- Instituto Superior T\'{e}cnico \\
% Av.\ Rovisco Pais, Torre Norte 7 \\
% 1049--001 Lisboa, \textsc{Portugal} \\
% Phone: +351 218418284 (POE), +351 218418296 (JG, JX), +351 218418053 (PO) \\
% Fax: +351 218418291 \\
% Email: \{poguz, jpg, jxavier, pjcro\}@isr.ist.utl.pt
\else
\begin{IEEEkeywords}
Simultaneous Localization and Tracking, Sensor Networks, Maximum Likelihood Estimation, Semidefinite Programming, Laplacian noise, Source Localization
\end{IEEEkeywords}
\fi
% \IEEEpeerreviewmaketitle

\section{Introduction}
\label{int}
\IEEEPARstart{T}{his} work addresses the problem of tracking a single target from distance-like measurements taken by nodes in a sensor network whose positions are not precisely known. The goal is to estimate the positions of all sensors and the target, given only partial or no {\em a priori} information on the spatial configuration of the network. The ability to track a target is a key component in several scenarios of wireless sensor networks, hence methods that avoid the need for careful calibration of sensor positions are practically relevant.

In \cite{Taylor2006,Funiak2006} SLAT is formulated in a Bayesian framework that resembles the related and well-studied problem of Simultaneous Localization and Mapping (SLAM) in robotics. The {\em a posteriori} probability density function of sensor/target positions and calibration parameters is recursively propagated in time as more target sightings become available. In \cite{Taylor2006} these observations are true range measurements obtained through a combination of transmitted acoustic and radio pulses. Alternatives to range include pseudorange and bearing information estimated from camera images \cite{Funiak2006}, or the (somewhat unreliable) Received Signal Strength (RSS) of radio transmissions \cite{Hero2008}. In \cite{Teng2009} the SLAT problem is also formulated in a Bayesian framework as a general state evolution model under a binary proximity model and solved in a decentralized way using binary sensor networks. Other SLAT-like approaches include localization (calibration) as presented in \cite{Moses2003}, where positions and orientations of unknown sources and sensors are centrally obtained via ML based on time-of-arrival and angle-of-arrival measurements. Several of the above references emphasize decentralized processing, in line with the local observation model.

When target dynamics are not accounted for, the SLAT problem may be thought of as a special case of Sensor Network Localization (SNL) with a limited set of intersensor measurements. In fact, EDM and related matrix completion methods based on Semidefinite Programming (SDP) have been adopted previously for static SNL (see \cite{Ye2006} and references therein). EDM completion for SLAT has been discussed in \cite{Hero2008}, although the authors ended up pursuing an alternative approximate completion approach based on a variant of Multidimensional Scaling (MDS). Underwater and underground scenarios with uncertainty in anchor positions are considered in \cite{Lui2009}, and edge-based SDP is proposed to reduce the computational complexity of SNL. In \cite{Weiss2008} static SNL is formulated as a problem of ML phase retrieval.

In addition to centralized SNL approaches such as \cite{Ye2006,Lui2009,Weiss2008}, enumerated above, a wealth of results is available on distributed approaches for scenarios where the existence of a central node is inconvenient, e.g., due to congested communications in its vicinity, or excessive vulnerability of the whole infrastructure to failure of that single node \cite{Hero2008,Luo2007,Vemula2009,So2009,Costa2006,Moura2009}. A two-step approach based on second-order cone programming relaxation with inaccurate anchor positions is introduced in \cite{Luo2007}. In \cite{So2009} a weighted least-squares algorithm with successive refinement provides both position estimates and their covariances in partially connected scenarios. A distributed weighted MDS method with majorization approximations is applied in \cite{Costa2006}. The cost function and majorization technique are similar to the ones used in this paper for ML iterative refinement under Gaussian noise, but initialization relies on prior estimates of sensor positions. A wholly different iterative approach for distributed SNL using barycentric coordinates and Cayley-Menger determinants is developed in \cite{Moura2009}. For $m$-dimensional Euclidean space a node reinterpolates its position based on estimates from a set of $m+1$ neighbors such that it lies in their convex hull.

This paper focuses on centralized SLAT based on plain ML estimation. A basic iterative optimization approach using Majorization Minimization (MM) is first developed for batch estimation, i.e., when all measurements are processed simultaneously. A time recursive method is then obtained by estimating each new target or sensor position as the corresponding range measurements become available, and then iteratively re-optimizing the expanded ML cost function. This continuation scheme only requires initialization of all target/sensor positions at the first time step, which is computationally less complex than doing so for each new sensed target or sensor position. EDM completion for batch estimation is proposed to initialize the iterative ML algorithm with little {\em a priori} knowledge of sensor/target positions, thus alleviating the problem of convergence to undesirable local extrema.

This approach was proposed in \cite{Ekim2009} for Gaussian noise, using cost functions for batch and time recursive initialization schemes that match \emph{squared} observations with \emph{squared} estimated ranges. These discrepancies are eliminated in the present paper, such that both initialization and ML refinement operate with cost functions that match plain (non-squared) ranges, leading to improved robustness under strong measurement noise \cite{Kay1993}. We use the Source Localization in Complex Plane (SLCP) approach proposed in \cite{Ekim2010} to obtain a cost function for initialization of incremental target/sensor position estimates that admits an accurate convex relaxation as an SDP. 

This work also develops modified versions of initialization and ML refinement for Laplacian noise, which models the presence of outliers in some practical ranging systems that adversely affect the performance of localization algorithms designed for Gaussian noise \cite{Taylor2006, Picard2009}. This is accomplished by replacing $\ell_{2}$ norms with $\ell_{1}$ norms in the optimization problems for initialization and refinement that were previously formulated for Gaussian noise, and then performing suitable manipulations to write these in a form that is amenable to general-purpose solvers. A source localization algorithm with $\ell_1$ norms (SL$\ell_1$) is also derived for incremental initialization of the SLAT scheme.

In \cite{Picard2009} $\ell_{1}$ norms are also used to handle outliers, but the proposed method is very different from the one developed here, as it relies on linear programming to identify the outliers, and then remove them from consideration when computing the source location. In our work all measurements are kept, as the modulus of range differences that appears in cost functions ensures that outlier terms do not overwhelm the remaining ones as long as the proportion of outliers remains small. Another approach for handling outliers is presented in \cite{Veen2009}, where the Huber cost function interpolates between $\ell_{1}$ and $\ell_{2}$ norms. This function is minimized via iterative majorization techniques with {\em a priori} information on sensor positions, where in each step a majorization subproblem is solved using Costa's algorithm.

The paper is organized as follows. In section \ref{Problem Formulation}, the SLAT problem is introduced. Section \ref{sec_gaussnoise} presents estimation methods for range measurements corrupted by Gaussian noise, namely, EDM initialization, iterative likelihood refinement by MM, and time-recursive updating through incremental estimation of target/sensor positions. Section \ref{sec_laplacenoise} develops similar methods for Laplacian noise. Section \ref{res} provides simulation results of batch and time recursive approaches. In addition, the performances of initialization techniques are compared with and without outliers. Finally, section \ref{conc} summarizes the main conclusions.

Throughout, both scalars and individual position vectors in (2D) space will be represented by lowercase letters. Matrices and vectors of concatenated coordinates will be denoted by boldface uppercase and lowercase letters, respectively. The superscript \emph{T} (\emph{H}) stands for the transpose (Hermitian) of the given real (complex) vector or matrix. Below, $\mathbf{I}_m$ is the $m \times m$ identity matrix and $\mathbf{1}_m$ is the vector of $m$ ones. For symmetric matrix $\mathbf{X}$, $\mathbf{X}\succeq 0$ means that $\mathbf{X}$ is positive semidefinite. 

%%%%%%%%%%%%%%%%%%%%%%%%%
\section{Problem Formulation}\label{Problem Formulation}
%%%%%%%%%%%%%%%%%%%%%%%%%

The network comprises sensors at unknown positions $\{x_{1},x_{2},\ldots,x_{n}\} \in \mathbb{R}^2$, a set of reference sensors (anchors) at known positions $\{a_{1},a_{2},\ldots,a_{l} \}\in \mathbb{R}^2$, and unknown target positions  $\{e_{1},e_{2},\ldots,e_{m} \}\in \mathbb{R}^2$. A central processing node has access to range measurements between target positions and all sensors/anchors, namely, 
\begin{align*}
d_{ij} & = \|x_{i}-e_{j}\| +w_{ij}, & d_{kj} & = \|a_{k}-e_{j}\| +w_{kj},
\end{align*} 
where  $w_{ij}$ and $w_{kj}$ denote noise terms. A practical system that provides such range measurements is used, e.g., in \cite{Taylor2006}.

\paragraph{SLAT Under Gaussian Noise} If disturbances are Gaussian, independent and identically distributed (i.i.d.), maximizing the likelihood for the full batch of observations is equivalent to minimizing the cost function
\begin{equation}\label{costfunc}
\Omega_G ({\bf x}) = \sum_{i,j}(\|x_{i}-e_{j}\|-d_{ij})^{2} + \sum_{k,j}(\|a_{k}-e_{j}\|-d_{kj})^{2}.
\end{equation}
The full set of unknown sensor and target positions is concatenated into column vector ${\bf x}\in \mathbb{R}^{2(n+m)}$, the argument of $\Omega_G $. The goal of our SLAT approach is to find the set of coordinates in $\mathbf{x}$ which minimizes \eqref{costfunc}.

\paragraph{SLAT Under Laplacian Noise} When the disturbances are Laplacian and i.i.d., thus heavier tailed than Gaussian, maximizing the likelihood amounts to minimizing the cost function
\begin{equation}\label{costfuncnew}
\Omega_L ({\bf x}) = \sum_{i,j}|\|x_{i}-e_{j}\|-d_{ij}| + \sum_{k,j}|\|a_{k}-e_{j}\|-d_{kj}|.
\end{equation}
When compared with \eqref{costfunc}, the absence of squares in the summation terms of \eqref{costfuncnew} renders the function less sensitive to outlier measurements $d_{ij}$ with large deviations from the true ranges.

Due to the nature of this problem the functions $\Omega_G $ and $\Omega_L$ are invariant to global rotation, translation and reflection in the absence of anchors. In order to remove ambiguities in the solutions, at least ~$l = 3$ non collinear anchors must be considered. As in many other ML problems, the functions $\Omega_G$ and $\Omega_L$ are in general nonconvex and multimodal, hence their (approximate) minimization proceeds in two phases, initialization and refinement. The former provides a suitable initial point through EDM completion, which tends to avoid convergence towards undesirable local minimizers of the ensuing iterative refinement algorithms based on MM or weighted-MM.

%%%%%%%%%%%%%%%%%%%%%%%%%
\section{SLAT under Gaussian Noise}\label{sec_gaussnoise}
%%%%%%%%%%%%%%%%%%%%%%%%%

This section develops algorithms for EDM initialization, MM refinement, and recursive estimation in  SLAT under the assumption that measurement noise is i.i.d.\ and Gaussian. Modifications for i.i.d.\ Laplacian noise are considered in Section \ref{sec_laplacenoise}. A basic formulation of EDM completion with squared distances is provided first, and will form the basis for the initialization methods described in Sections \ref{EDM_pd} and \ref{EDM_l1}.

\subsection{EDM with Squared Distances}\label{EDM_sr}
The basic EDM completion problem, described below, operates on squared ranges \cite{Boyd2004,Dattorro2005}. Even though it is not matched to the likelihood function \eqref{costfunc}, it will be useful for benchmarking in Section \ref{res}, as its performance is representative of other popular SNL methods \cite{Ye2006,Ding2006} and the SLAT approach of \cite{Ekim2009}.

A partial pre-distance matrix ${\bf D}$ is a matrix with zero diagonal entries and with certain elements fixed to given nonnegative values which are the squared observed distances, $D_{ij}=d_{ij}^{2}$. The remaining elements are considered free. The nearest EDM problem is to find an EDM {\bf E} that is nearest in the least-squares sense to matrix ${\bf D}$, when the free variables are not considered and the elements of {\bf E} satisfy $E_{ij} = \|y_i-y_j\|^2$ for a set of points ${y_{i}}$. The geometry and properties of EDM (a convex cone) have been extensively studied in the literature \cite{Boyd2004,Dattorro2005}. The nearest EDM problem in 2D space is formulated as

\begin{equation}\label{approx1EDM}
\begin{array}{cl}
\mbox{minimize} &\sum_ {i,j \in \cal{O}}(E_{ij}-d_{ij}^2)^{2}\\
{\bf E}\\
\mbox{subject to} & {\bf E}\in \edmc \\
&\mathbf{E}(\cal{A}) = \mathbf{A} \\
&\mbox{rank}({\bf JEJ} ) = 2,
\end{array}
\end{equation}
where
\begin{equation*}
{\bf J}=\bigl( {\bf I}_{\rho}-\frac{1}{\rho} {\bf 1}_{\rho} {{\bf 1}_{\rho}}^T\bigr ), \quad \rho = m+n+l,
\end{equation*}
is a centering operator which subtracts the mean of a vector from each of its components. In \eqref{approx1EDM}, $\cal{O}$ is the index set for which range measurements are available. The constraint $\mathbf{E}(\cal{A}) = \mathbf{A}$, where $\cal{A}$ is the index set of anchor/anchor distances and $\mathbf{A}_{ij} = \|a_{i}-a_{j}\|^2$ is the corresponding EDM submatrix, enforces the known \emph{a priori} spatial information. Matrix {\bf E} belongs to the EDM cone $\edmc$ if it satisfies the properties 
\begin{align}
E_{ii} & = 0, & E_{ij} & \geq 0, & {\bf -JEJ} & \succeq 0.
\end{align}
The rank constraint in (\ref{approx1EDM}) ensures that the solution is compatible with a constellation of sensor/anchor/target points in $\mathbb{R}^2$. Extraction of the set ${y_{i}}$ from {\bf E} is described below. Problem (\ref{approx1EDM}) is also known as the penalty function approximation \cite{Boyd2004} due to the form of the cost function $\varphi_1({\bf E}) = \sum_{i,j}(E_{ij}-d_{ij}^2)^{2}$.
Expressing (\ref{approx1EDM}) in terms of full matrices and dropping the rank constraint, a compact relaxed SDP formulation is obtained as 
\begin{equation}\label{approx3EDM}
\begin{array}{cl}
\mbox{minimize}&\|{\bf W}\odot{\bf (E-D)}\|_{F}^{2}\\
{\bf E}\\
\mbox{subject to}&{\bf E}\in \edmc, \; \mathbf{E}(\cal{A}) = \mathbf{A},
\end{array}
\end{equation}
where ${\bf W}$ is a mask matrix with zeros in the entries corresponding to free elements of $D_{ij} = d_{ij}^2$, and ones elsewhere. Combined with the Hadamard product $\odot$, and Frobenius norm $\|.\|_{F}$, this replaces the summation in \eqref{approx1EDM} over the observed index set $\cal{O}$.
From here on, we will call this method EDM with squared ranges (EDM-SR). 

\subsection{SLAT Initialization: EDM with Plain Distances}\label{EDM_pd}
Instead of trying to match squared distances, we can apply EDM completion to plain distances as \begin{equation}\label{approx1met2}
\begin{array}{cl}
\mbox{minimize}&\sum_ {i, j} (\sqrt{E_{ij}}-d_{ij})^{2}\\
{\bf E}\\
\mbox{subject to}&{\bf E}\in \edmc , \; \mathbf{E}(\cal{A}) = \mathbf{A}\\
&\mbox{rank}({\bf JEJ}) = 2.
\end{array}
\end{equation}
For this method the penalty function is $\varphi_2({\bf E}) = \sum_{i,j}(\sqrt{E_{ij}}-d_{ij})^{2}$, which more closely matches the terms in the likelihood function (\ref{costfunc}), and (\ref{approx1met2}) is thus expected to inherit some of the robustness properties of ML estimation. Expanding the objective function in (\ref{approx1met2}) results in
\begin{equation}\label{approx2met2}
\begin{array}{cl}
\mbox{minimize}&\sum_ {i, j} (E_{ij}-2\sqrt{E_{ij}}d_{ij}+d_{ij}^2)\\
{\bf E}\\
\mbox{subject to}&{\bf E}\in \edmc , \; \mathbf{E}(\cal{A}) = \mathbf{A}\\
&\mbox{rank}({\bf JEJ}) = 2.
\end{array}
\end{equation}
A relaxed SDP is obtained by introducing an epigraph-like variable ${\bf T}$ and dropping the rank constraint
\begin{equation}\label{approx4met2}
\begin{array}{cl}
\mbox{minimize}&\sum_ {i, j} (E_{ij}-2T_{ij}d_{ij})\\
{\bf E,T}\\
\mbox{subject to}& T_{ij}^2 \leq E_{ij}\\
&{\bf E}\in \edmc , \; \mathbf{E}(\cal{A}) = \mathbf{A}.
\end{array}
\end{equation}
From here on, we will call this method EDM with plain ranges (EDM-R).

Remark that the solutions of the initialization techniques described here and in Sections \ref{EDM_sr} and \ref{EDM_l1} are distance matrices. Detailed explanations of how to estimate the spatial coordinates of the sensors and target positions from EDM and the usage of anchors are given in \cite{Ekim2009}. The basic idea is to use a linear transformation to obtain the Gram matrix ${\bf Y}^{T}{\bf Y}$ from the EDM matrix {\bf E}, from which spatial coordinates {\bf Y} are extracted by singular value decomposition (SVD) up to a unitary matrix. The anchors are then used to estimate the residual unitary matrix by solving a Procrustes problem. As discussed in \cite{Ekim2009}, observation noise can significantly disrupt the estimated sensor/target coordinates through EDM completion and rank truncation, and it was found that much more accurate results are obtained by using those as a starting point for likelihood maximization. MM algorithms are proposed next for iterative likelihood maximization.

\subsection{SLAT Iterative Refinement: Majorization-Minimization}\label{maj-min}

The key idea of MM is to find, at a certain point ${\bf x}^{t}$, a simpler function that has the same function value at ${\bf x}^{t}$ and anywhere else is larger than or equal to the objective function to be minimized. Such a function is called a majorization function. By minimizing the majorization function we obtain the next point of the algorithm, which also decreases the cost function \cite{Hunter2004}. The detailed derivation of MM is given in \cite{Ekim2009}, but it is summarized below for the reader's convenience. Define two convex functions as
\begin{align} \label{eq_dist}
f_{ij}({\bf x}) & =\|x_{i}-e_{j}\|, & g_{kj}({\bf x}) & = \|a_{k}-e_{j}\|.
\end{align}
Expanding $f$ and $g$ in \eqref{costfunc} and using first-order conditions on convexity \cite{Boyd2004}, 
\begin{equation}\label{expfunc}
\begin{split} 
\Omega_G({\bf x}) \leq &\sum_{i,j}\left( f_{ij}^{2}({\bf x}) - 2d_{ij}\bigl( f_{ij}({\bf x}^{t}) + \langle \nabla f_{ij}({\bf x}^{t}), ({\bf x}-{\bf x}^{t}) \rangle \bigr) + d_{ij}^2 \right) \\
+& \sum_{k,j}\left( g_{kj}^{2}({\bf x}) - 2d_{kj}\bigl( g_{kj}({\bf x}^{t})+ \langle \nabla g_{kj}({\bf x}^t), ({\bf x}-{\bf x}^{t})\rangle  \bigr) + d_{kj}^2 \right),
\end{split}
\end{equation}
where $\langle u,v \rangle = u^Tv$, we get the proposed majorization function on the right side of (\ref{expfunc}), which is quadratic in ${\bf x}$ and easily minimized. Hence the MM iteration
\begin{equation}
{\bf x}^{t+1} = \arg\min_{{\bf x}} \sum_{i,j}\left( f_{ij}^{2}({\bf x})-2d_{ij} \langle \nabla f_{ij}({\bf x}^{t}), {\bf x}\rangle  \right)
+\sum_{k,j}\left( g_{kj}^{2}({\bf x})-2d_{kj} \langle \nabla g_{kj}({\bf x}^{t}), {\bf x}\rangle \right)
\end{equation}
turns out to be obtained as the solution of a linear system of equations. 

\subsection{Time-Recursive Position Estimation: SLCP}\label{timerec}

Suppose that a batch of observations has been processed, and a new target position is to be estimated.  We could repeat the previous approach by redefining the new batch as the old one concatenated with the new set of observations. However, this is computationally expensive due to the EDM step. Also, previous estimated positions would be ignored. To alleviate the load we propose a simple methodology to obtain a good initial point which avoids the EDM step. It consists of fixing the previous positions at their estimated values and only estimating the new target position. More precisely, we minimize
\begin{equation}\label{sqcostfunc}
\Psi_G (y) = \sum_{i = 1}^{n+l}(\|b_{i}-y\|-d_{i})^{2},
\end{equation}
where $y$ is the new target position, $b_i$ denotes the position of a sensor or anchor, and $d_i$ is the corresponding range measurement. We propose SLCP \cite{Ekim2010} to minimize \eqref{sqcostfunc}, and briefly summarize the method below.

One can readily see that each term in \eqref{sqcostfunc} quantifies the distance between the target location and a circle, centered at an anchor or sensor, with radius equal to the measured range, e.g., $\bigl|\|b_{i}-y\|-d_{i}\bigr|=\|y-y_{i}\|$, where $y_{i}$ is located at the intersection of the line connecting $b_{i}$ and $y$ with the circle $\{z:\|z-b_{i}\|=d_{i}\}$. An equivalent formulation is therefore 
\begin{equation}\label{sourcecostfunction1}
\begin{array}{cl}
\mbox{minimize}&\sum_{i = 1}^{n+l} \|y-y_{i}\|^{2}\\
y,y_i&\\
\mbox{subject to}&\|b_{i} - y_i\|=d_{i} ~~i = 1, \ldots, n+l.
\end{array}
\end{equation}
SLCP relies on a complex formulation of \eqref{sourcecostfunction1} for the 2D case, and the problem is then manipulated and relaxed to an SDP of the form
\begin{equation}\label{sourcecostfunction11}
\begin{array}{cl}
\mbox{maximize}&t+\frac{1}{n+l}{\bf r}^T \bm{\Phi}{\bf r}\\
\bm{\Phi},t\\
\mbox{subject to}&\bm{\Phi} \succeq 0, \; \; \phi_{ii}=1\\
&4{\bf c}^H\bm{\Phi}{\bf c}\geq t^2,
\end{array}
\end{equation}
where $\mathbf{r}$ and $\mathbf{c}$ are constant vectors built from sensor/anchor complex coordinates and measured radii. The solution of the optimization problem (\ref{sourcecostfunction11}) is a positive semidefinite matrix, hopefully with near-1 rank. Afterwards, the target coordinates are estimated by SVD of $\bm{\Phi}$ as described in \cite{Ekim2009}.  

After an optimal target position is obtained, we return to the cost function (\ref{costfunc}) or  (\ref{costfuncnew}) and iteratively refine all the estimates. This incremental or time recursive procedure can be applied to either new targets or sensors.

In a previous paper \cite{Ekim2009}, the source localization method derived in \cite{Stoica2008}, termed Squared Range Least Squares (SR-LS), is proposed as the time recursive method. Note that in \cite{Stoica2008} \emph{squared} distances are matched, leading to a Trust Region optimization problem. However, as demonstrated in \cite{Ekim2010}, SLCP is a more accurate source localization method and its cost function \eqref{sqcostfunc} is better matched to the likelihood function \eqref{costfunc}. This makes it more convenient for initialization of iterative refinement algorithms, which therefore require fewer iterations to converge and/or are less likely to be trapped in undesirable local extrema.

%%%%%%%%%%%%%%%%%%%%%%%%%
\section{SLAT under Laplacian Noise}\label{sec_laplacenoise}
%%%%%%%%%%%%%%%%%%%%%%%%%

\subsection{SLAT Initialization: EDM with Ranges and $\ell_1$-norm}\label{EDM_l1}

\begin{figure}%[h]
\begin{center}
\includegraphics[width=0.6\textwidth]{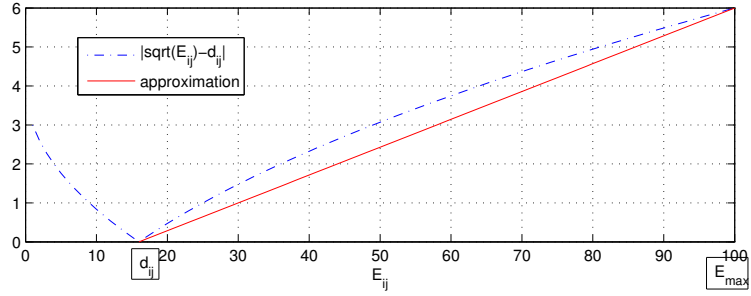}
\end{center}
\caption{The value of $\abs{\sqrt{E_{ij}}-d_{ij}}$ vs $E_{ij}$, and the linear approximation of the concave part.}
\label{approxfun}
\end{figure}

Among the penalty function approximation methods $\ell_1$-norm approximation is known to be robust to outliers \cite{Boyd2004}. Therefore, the penalty function of the third SLAT initialization method is chosen as $\varphi_3({\bf E}) = \sum_{i,j}\abs{\sqrt{E_{ij}}-d_{ij}}$, and the associated optimization problem becomes
\begin{equation}\label{approx13rdmet}
\begin{array}{cl}
\mbox{minimize}&\sum_ {i, j}\abs{\sqrt{E_{ij}}-d_{ij}}\\
{\bf E}\\
\mbox{subject to}&{\bf E}\in \edmc, \; \mathbf{E}(\cal{A}) = \mathbf{A}\\
&\mbox{rank}({\bf JEJ}) = 2.
\end{array}
\end{equation}
This problem is not convex, as the objective function $\abs{\sqrt{E_{ij}}-d_{ij}}$ is convex when $\sqrt{E_{ij}}-d_{ij}<0$, but concave for $\sqrt{E_{ij}}-d_{ij}>0$. To obtain a convex approximation the objective function is replaced by a linear approximation
\begin{align}
&  a_{ij}E_{ij}+b_{ij}, & a_{ij} & = \frac{1}{\sqrt{E_{\text{max}}}+d_{ij}}, & b_{ij} = -\frac{d_{ij}^2}{\sqrt{E_{\text{max}}}+d_{ij}}
\end{align}
 in part of the domain where it is concave, as shown in Figure \ref{approxfun}.
The two functions coincide for $E_{ij} = d_{ij}^{2}$ and $E_{ij} = E_{\text{max}}$, where the constant $E_{\text{max}}$ is a practical upper bound on (squared) range measurements. Thus, we replace $\abs{\sqrt{E_{ij}}-d_{ij}}$ by its convex envelope $\max\{d_{ij}-\sqrt{E_{ij}},a_{ij}E_{ij}+b_{ij}\}$ and use the epigraph variable ${\bf T}$ to obtain
\begin{equation}\label{approx23rdmet}
\begin{array}{cl}
\mbox{minimize}&\sum_ {i, j}  T_{ij}\\
{\bf E,T}\\
\mbox{subject to}&\max\{d_{ij}-\sqrt{E_{ij}},a_{ij}E_{ij}+b_{ij}\}\leq T_{ij}\\
&{\bf E}\in \edmc, \; \mathbf{E}(\cal{A}) = \mathbf{A}\\
&\mbox{rank}({\bf JEJ}) = 2.
\end{array}
\end{equation}
A relaxation of (\ref{approx23rdmet}) after dropping the rank constraint is
\begin{equation}\label{approx43rdmet}
\begin{array}{cl}
\mbox{minimize}&\sum_ {i, j}  T_{ij}\\
{\bf E,T}\\
\mbox{subject to} &(d_{ij}-T_{ij})^2\leq E_{ij}\\
&a_{ij}E_{ij}+b_{ij}\leq T_{ij}\\
&{\bf E}\in \edmc, \; \mathbf{E}(\cal{A}) = \mathbf{A}.
\end{array}
\end{equation}
Note that the first constraint in \eqref{approx43rdmet} is not equivalent to $d_{ij}-\sqrt{E_{ij}} \leq T_{ij}$, but rather to $-\sqrt{E_{ij}} \leq d_{ij} - T_{ij} \leq \sqrt{E_{ij}}$, which amounts to intersecting the original epigraph with the parabolic hypograph $d_{ij} + \sqrt{E_{ij}} \geq T_{ij}$. This preserves the convexity of the feasible set and does not change its lower boundary for $E_{ij} \in [0, \, E_{\text{max}}]$, where the optimal point will be found. The constraint can now be readily expressed in standard form without introducing additional variables, e.g., as an LMI or a second-order cone constraint \cite{Lobo1998}
\begin{align}
  \begin{bmatrix} 1 & d_{ij}-T_{ij} \\ d_{ij}-T_{ij} & E_{ij} \end{bmatrix} & \succeq 0 & \text{or} &&
  \begin{Vmatrix} 2(d_{ij}-T_{ij}) \\ E_{ij}-1 \end{Vmatrix} & \leq E_{ij}+1.
\end{align}
This technique will be called EDM with ranges and $\ell_1$-norm (EDM-R-$\ell_1$).

\subsection{SLAT Iterative Refinement: Weighted Majorization Minimization}\label{weighted maj_min}

Robustness to outliers in the cost function \eqref{costfuncnew} for Laplacian noise is gained at the expense of differentiability. To circumvent that shortcoming we resort to the well-known re-weighted least squares approach \cite{Moon2000}, which replaces the minimization of \eqref{costfuncnew} with a sequence of minimizations of smooth approximation functions that converge to $\Omega_{L} ({\bf x})$. Specifically, \eqref{costfuncnew} is first written as
\begin{equation}\label{costfuncweighted}
\Omega_{L} ({\bf x}) = \sum_{i,j}u_{ij}(\|x_{i}-e_{j}\|-d_{ij})^2 + \sum_{k,j}v_{kj}(\|a_{k}-e_{j}\|-d_{kj})^2,
\end{equation}
with
\begin{align*}
u_{ij}  & = \frac{1}{|\|x_{i}-e_{j}\|-d_{ij}|}, & v_{kj} & = \frac{1}{|\|a_{k}-e_{j}\|-d_{kj}|}.
\end{align*}
At time $t$ the function to be minimized becomes $\Omega_{L}^{t}(\mathbf{x})$, which has the same form of \eqref{costfuncweighted} but the \emph{functions} $u_{ij}$, $v_{kj}$ above are now replaced by \emph{constants} based on the estimated positions after the previous iteration
\begin{align} \label{eq_wmmweights}
u_{ij}^{t}  & = \frac{1}{|\|x_{i}^{t}-e_{j}^{t}\|-d_{ij}|}, & v_{kj}^{t} & = \frac{1}{|\|a_{k}-e_{j}^{t}\|-d_{kj}|}.
\end{align}
An inner optimization loop could now be used to minimize $\Omega_{L}^{t}(\mathbf{x})$ for every $t$ but, as shown in Appendix \ref{sec_wmmconv}, a single iteration suffices to ensure convergence. With fixed $u_{ij}^{t}$, $v_{kj}^{t}$ the same majorization technique of Section \ref{maj-min} then yields the weighted-MM iteration
\begin{equation}\label{itwMM}
{\bf x}^{t+1} = \arg\min_{{\bf x}}  \sum_{i,j} u_{ij}^{t}\left( f_{ij}^{2}({\bf x})-2d_{ij} \langle \nabla  f_{ij}({\bf x}^{t}), {\bf x} \rangle  \right)
+\sum_{k,j} v_{kj}^{t}\left( g_{kj}^{2}({\bf x})-2d_{kj} \langle \nabla  g_{kj}({\bf x}^{t}), {\bf x} \rangle \right).
\end{equation}
In practice the weights $u_{ij}^t$ and $v_{ij}^t$ must be modified to avoid the possibility of division by zero \cite{Rodriguez2009}, which in our case is achieved by saturating them at $10^{5}$ when computing \eqref{itwMM}.

\subsection{Time-Recursive Position Estimation: SL$\ell_1$}\label{sec:timerec_l1}

The ML source localization problem under Laplacian noise is equivalent to
\begin{equation}\label{sourcecostfuncl1}
\begin{array}{cl}\mbox{minimize}& \Psi_{L}(y) = \sum_{i=1}^{n+l} |\|y-b_i\|-d_{i}|\\
y
\end{array}
\end{equation}
or
\begin{equation}\label{sourcecostfuncl1.2}
\begin{array}{cl}
\mbox{minimize}& (\sum_{i=1}^{n+l} |\|y-b_i\|-d_{i}|)^2,\\
y
\end{array}
\end{equation}
where $y$, $b_i$ and $d_i$ are defined in section \ref{timerec}. We use ideas from \cite{Nemirovsky2001} to express the minimization of $\Psi_{L}^2$ as a weighted sum of squares.

\begin{thm} \label{thm:1}
The following problem is equivalent to \eqref{sourcecostfuncl1.2} 
\begin{equation}\label{sourcecostfuncl1.3}
\begin{array}{ccl}
\mathrm{minimize} & \mathrm{minimize} & \sum_{i=1}^{n+l} \frac{(\|y-b_i\|-d_{i})^2}{\lambda_i}\\
y & \bm{\lambda} \in \mathbb{R}^{n+l}\\
& \mathrm{subject~to} & \lambda_i>0,~ \mathbf{1}^T\bm{\lambda} = 1.
\end{array}
\end{equation}
\end{thm}
A proof is given in Appendix \ref{sec_sll1proofs}. As claimed in section \ref{timerec} the difference between the true range and observed range is actually equivalent to the distance between the source position and the point on the circle with center $b_i$ and radius $d_i$. An equivalent formulation is therefore
\begin{equation}\label{sourcecostfuncl1.4}
\begin{array}{cl}
\mbox{minimize}&\sum_{i=1}^{n+l} \frac{\|y-y_i\|^{2}}{\lambda_i}\\
y,y_i,\bm{\lambda}\\
\text{subject to} & \|y_i-b_i\|=d_i \\
& \lambda_i>0,~ \mathbf{1}^T\bm{\lambda}=1. 
\end{array}
\end{equation}
If we fix the $y_{i}$ and $\lambda_i$, the solution of \eqref{sourcecostfuncl1.4} with respect to $y$ is an unconstrained optimization problem whose solution is readily obtained by invoking the optimality conditions 
\begin{equation} \label{sl1_targetcoords}
\sum_{i=1}^{n+l}\frac{(y-y_i)}{\lambda_i}=0 \; \Rightarrow \; y^*=\frac{\sum_{i=1}^{n+l}\frac{y_i}{\lambda_i} }{\sum_{i=1}^{n+l}\frac{1}{\lambda_i}}.
\end{equation}
Geometrically, the first constraint of \eqref{sourcecostfuncl1.4} defines circle equations, which can be compactly described in the complex plane as $y_i = b_i+d_ie^{j\phi_i}$. We collect these into a vector $\mathbf{y = b + Ru}$, where $\mathbf{b}= \begin{bmatrix} b_1 & \ldots & b_{n+l}\end{bmatrix}^T \in \mathbb{C}^{n+l}$, $\mathbf{R}= \text{diag}(d_1, \ldots, d_{n+l})\in \mathbb{R}^{(n+l) \times (n+l)}$, and $\mathbf{u} = \begin{bmatrix} e^{j\phi_1} & \ldots & e^{j\phi_{n+l}} \end{bmatrix}^T \in \mathbb{C}^{n+l}$. Using the optimal $y$, we get
\begin{equation}\label{sourcecostfuncl1.5}
\begin{array}{cl}
\mbox{minimize}&{\bf y}^H\bm{\Pi}{\bf y} = (\mathbf{b+Ru})^H\bm{\Pi}(\mathbf{b+Ru})\\
\bm{\lambda}, \mathbf{u}\\
\text{subject to} & \lambda_i>0,~ \mathbf{1}^T\bm{\lambda}=1\\
& \abs{u_i} = 1,
\end{array}
\end{equation}
where
\begin{equation}\label{sll1_projector}
\begin{split}
\bm{\Pi} & = \left[\begin{array}{ccc}
\frac{1}{\lambda_1}&0&0\\ 
0& \ddots &0\\
0& 0 & \frac{1}{\lambda_{n+l}}\end{array}\right]
-\frac{1}{\sum_{i=1}^{n+l}\frac{1}{\lambda_i}}
\left [
\begin{array}{c}
\frac{1}{\lambda_1}\\
\vdots\\
\frac{1}{\lambda_{n+l}} \end{array}\right]
\left [
\begin{array}{ccc}
\frac{1}{\lambda_1}&\ldots&\frac{1}{\lambda_{n+l}} \end{array}\right]\\
& = \bm{\Lambda}^{-1} - \bm{\Lambda}^{-1}\mathbf{1} (\mathbf{1}^{T}\bm{\Lambda}^{-1}\mathbf{1})^{-1}\mathbf{1}^{T}\bm{\Lambda}^{-1},
\end{split}
\end{equation}
with $\bm{\Lambda} = \text{diag}(\lambda_1,\ldots,\lambda_{n+l})$.

Matrix $\bm{\Pi}$ resembles an orthogonal projector. Using the matrix inversion lemma\footnote{$(A+BCD)^{-1} = A^{-1}-A^{-1}B(DA^{-1}B+C^{-1})^{-1}DA^{-1}$.} it is seen to be the limiting case $\bm{\Pi} = \lim_{\sigma \to \infty}(\bm{\Lambda} + \sigma \mathbf{1} \mathbf{1}^{T})^{-1}$ and thus positive semidefinite. This format is more amenable to analytic manipulations in optimization problems and will be used throughout. The parameter $\sigma$ is taken as a sufficiently large constant (see Appendix \ref{sec_sll1proofs}), although it could also be regarded as an additional optimization variable to ensure adequate approximation accuracy.

We now introduce an epigraph variable $t \in \mathbb{R}$ in \eqref{sourcecostfuncl1.5}, i.e., we minimize over $t$ and add the constraint $t -  (\mathbf{b+Ru})^H\bm{\Pi}(\mathbf{b+Ru}) \geq 0$.  Applying Schur complements the constraint may be successively written as
\begin{align}
  \begin{bmatrix}
    t & (\mathbf{b+Ru})^H \\
    \mathbf{b+Ru} & \bm{\Pi}^{-1}
  \end{bmatrix} & \succeq 0 & \Leftrightarrow & & \bm{\Pi}^{-1} -\frac{(\mathbf{b+Ru})(\mathbf{b+Ru})^H}{t} & \succeq 0.
\end{align}
The formulation becomes
\begin{equation}\label{sourcecostfuncl1.6}
\begin{array}{cl}
\text{minimize} & t\\
\bm{\lambda}, \mathbf{u}, t\\
\text{subject to} & \lambda_i>0,~ \mathbf{1}^T\bm{\lambda}=1\\
& \abs{u_i} = 1\\
& t \bm{\Lambda} + t \sigma \mathbf{1} \mathbf{1}^{T} \succeq (\mathbf{b+Ru})(\mathbf{b+Ru})^H.
\end{array}
\end{equation}
Finally, we define $\mathbf{B} = [\mathbf{b}~~ \mathbf{R}]$, $\mathbf{v}^H = [1 ~~\mathbf{u}^H]$, $\mathbf{V} = \mathbf{vv}^H$, and ignore the rank-1 constraint on the new variable $\mathbf{V}$ to obtain the relaxed SDP
\begin{equation}\label{sourcecostfuncl1.7}
\begin{array}{cl}
\text{minimize} & t\\
\bm{\beta}, \mathbf{V}, t\\
\text{subject to} & \beta_i>0,~ \mathbf{1}^T\bm{\beta}=t\\
& V_{ii} = 1, \mathbf{V} \succeq 0\\
& \text{diag}(\bm{\beta}) + t \sigma \mathbf{1} \mathbf{1}^{T} \succeq \mathbf{BVB}^H.
\end{array}
\end{equation}
The solution of the optimization problem \eqref{sourcecostfuncl1.7} includes the positive semidefinite matrix $\mathbf{V}$ from whose first row or column $\mathbf{u}$ can be extracted directly\footnote{Alternatively, $\mathbf{u}$ can be obtained by rank 1 factorization of the lower right submatrix of $\mathbf{V}$ corresponding to $\mathbf{uu}^H$, as in \cite{Ding2006, Ye2006}.}. to obtain $\mathbf{y = b+Ru}$ and the target coordinates from \eqref{sl1_targetcoords}.

As in section \ref{timerec}, after an optimal target position is obtained we return to the cost function \eqref{costfuncnew} and iteratively refine all the estimates. It is demonstrated in simulation in section \ref{res} that SL$\ell_1$ is more robust to outliers than the SLCP algorithm of section \ref{timerec}, as its cost function \eqref{sourcecostfuncl1} is better matched to the likelihood function \eqref{costfuncnew}.   

%%%%%%%%%%%%%%%%%%%%%%%%%%%%%%%%%%%%%%%%%%%%%%%%%%%%%%%%%%%
\section{Numerical Results}\label{res}
\begin{figure}%[h]
\centering
\subfloat[Without outliers\label{compinit:nooutliers}]{\includegraphics[width=5.5in, height=3in]{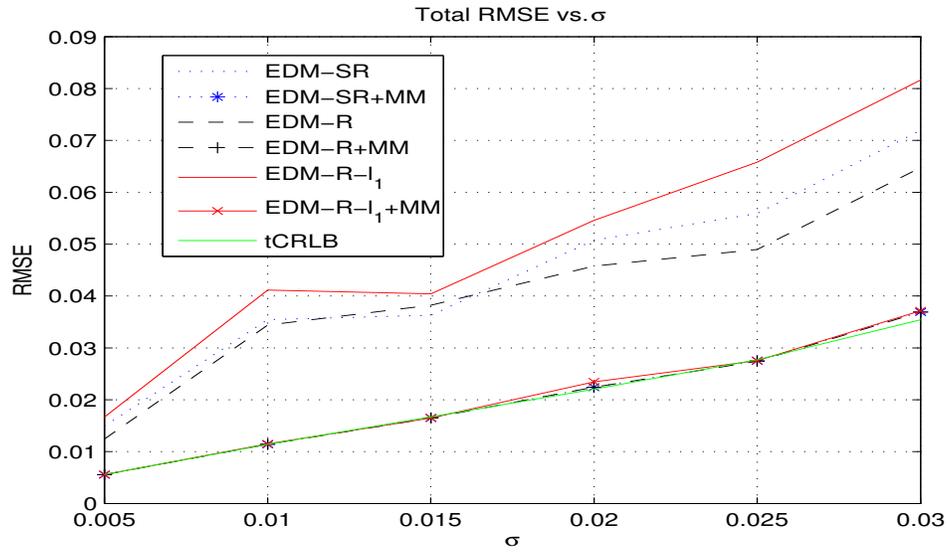}}\\
\subfloat[With outliers, $\sigma_{\text{gaussian}} = 0.01$\label{compinit:outliers}]{\includegraphics[width=5.5in, height=3in]{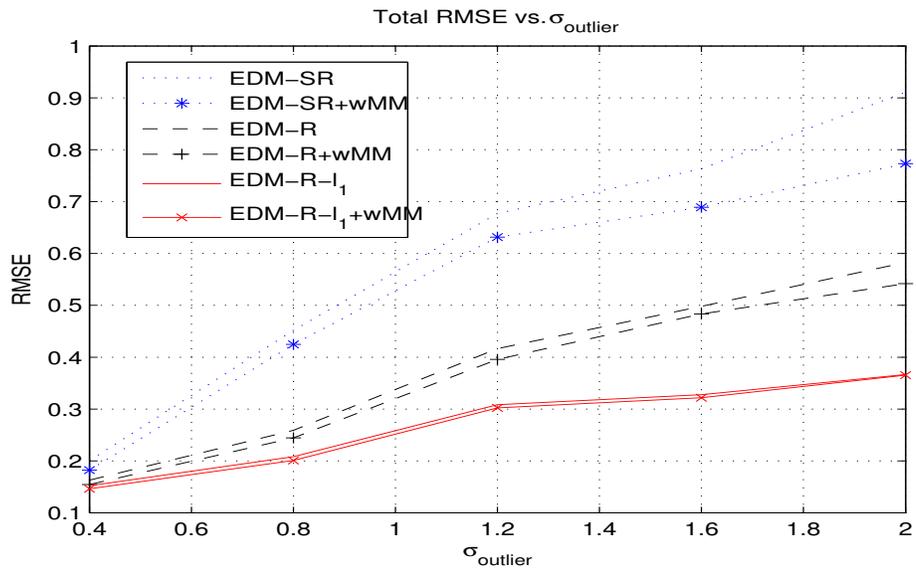}}
\caption{Comparison of batch initialization and refinement methods for SLAT.}
\label{compinit}
\end{figure}

%%%%%%%%%%%%%%%%%%%%%%%%%%%%%%%%%%%%%%%%%%%%%%%%%%%%%%%%%%%%%%%%%%
\paragraph*{\textbf{Example 1 [Comparison of Batch Initialization Methods]}} To investigate the accuracy of the methods, we set a physical scenario containing four anchors, five unknown sensors, and six target positions in a $[0,2]\times [0,2]$ area. Range measurements are corrupted by additive spatio-temporally white noise with standard deviation $\sigma_{\text{gaussian}} \in [0.005, 0.03]$. This noisy observation model may lead to near-zero or negative range measurements, in which case we follow normal practice \cite{Ye2006} and set them equal to a small positive constant ($10^{-5}$ in our simulations). With the chosen noise variances this occurs sufficiently seldom (up to 0.04\% of measurements) for its impact on estimation accuracy to be unimportant. Several algorithms are tested (EDM-SR, EDM-R, EDM-R-$l_{1}$, MM initialized by EDM-SR (EDM-SR+MM), MM initialized by EDM-R (EDM-R+MM) and MM initialized by EDM-R-$l_{1}$ (EDM-R-$l_{1}$+MM)), and their performances compared according to the total root mean square error (RMSE) 
\begin{equation}
\sqrt{\frac{1}{K}\frac{1}{n+m}\sum_{k=1}^K\sum_{i=1}^{n+m}\|x_i-\hat{x}_i^k\|^2},
\end{equation}
where $\hat{x}_i^k$ denotes the $i$-th estimated sensor or target position in the $k$-th Monte Carlo run for the specific noise realization. In each of $K=150$ Monte Carlo runs, a random network is generated according to the physical scenario described above. To assess the fundamental hardness of position estimation, error plots for Gaussian noise also show the total Cram\'{e}r-Rao Lower Bound (CRLB), calculated as
\begin{equation}\label{crbound}
\sqrt{\frac{1}{K}\frac{1}{n+m}\sum_{k=1}^K \text{trace}(\text{CRLB}_k)}
\end{equation}
for each noise variance, where $\text{CRLB}_k$ denotes the $k$-th diagonal element of the matrix lower bound. The CRLB for anchored and anchor-free localization using ranging information has been studied in \cite{Patwari2003,Chang2006,Jia2008} for different variance models of range estimation noise. For convenience, the CRLB for our SLAT problem under Gaussian noise is rederived in Appendix \ref{sec_crlb} in terms of the notation adopted in this paper. We do not prove unbiasedness of our estimators, a mathematically challenging endeavor that would be required to fully justify benchmarking against the CRLB. In our experimental results, however, we found no clear evidence of bias for small noise levels, where convergence to undesirable extrema of the cost functions is avoided. Figure~\ref {compinit:nooutliers} shows that plain EDM-R has better accuracy than EDM-SR and EDM-R-$l_{1}$, although the performance gap closes after iterative refinement by MM. Moreover, MM initialized by the various methods nearly touches the CRLB except when the noise variance is large. 

To compare the total RMSE of the algorithms in the presence of outliers, modified range measurements are created according to  a ``selective Gaussian'' model $d_{i} = \|\cdot\|+w_{i}+|\epsilon_i|$, where $\epsilon_i$ is a white Gaussian noise term with standard deviation $\sigma_{\text{outlier}}\in [0.4,2]$. The disturbance $\epsilon_i$ randomly affects only two range measurements, whereas $w_{i}$ with $\sigma_{\text{gaussian}} = 0.01$ is present in all observations. This outlier generating model deviates from the earlier Laplacian asssumption, but it is arguably representative of observed range measurements in practical systems \cite{Weiss2008}. Numerical results under a pure Laplacian model will be presented in Examples 3 and 4. In the presence of high noise and/or outliers, Fig.~\ref {compinit:outliers} shows that weighted-MM refinement does not close the performance gap between EDM-R-$l_{1}$, EDM-R and EDM-SR initialization because in the latter cases the algorithms converge more often to local minima, thus producing a larger total RMSE.
%%%%%%%%%%%%%%%%%%%%%%%

\begin{figure}%[h]
\begin{center}
\subfloat[\label{mean_var_wholeview}Whole constellation]{\includegraphics[width=5.5in, height=3in]{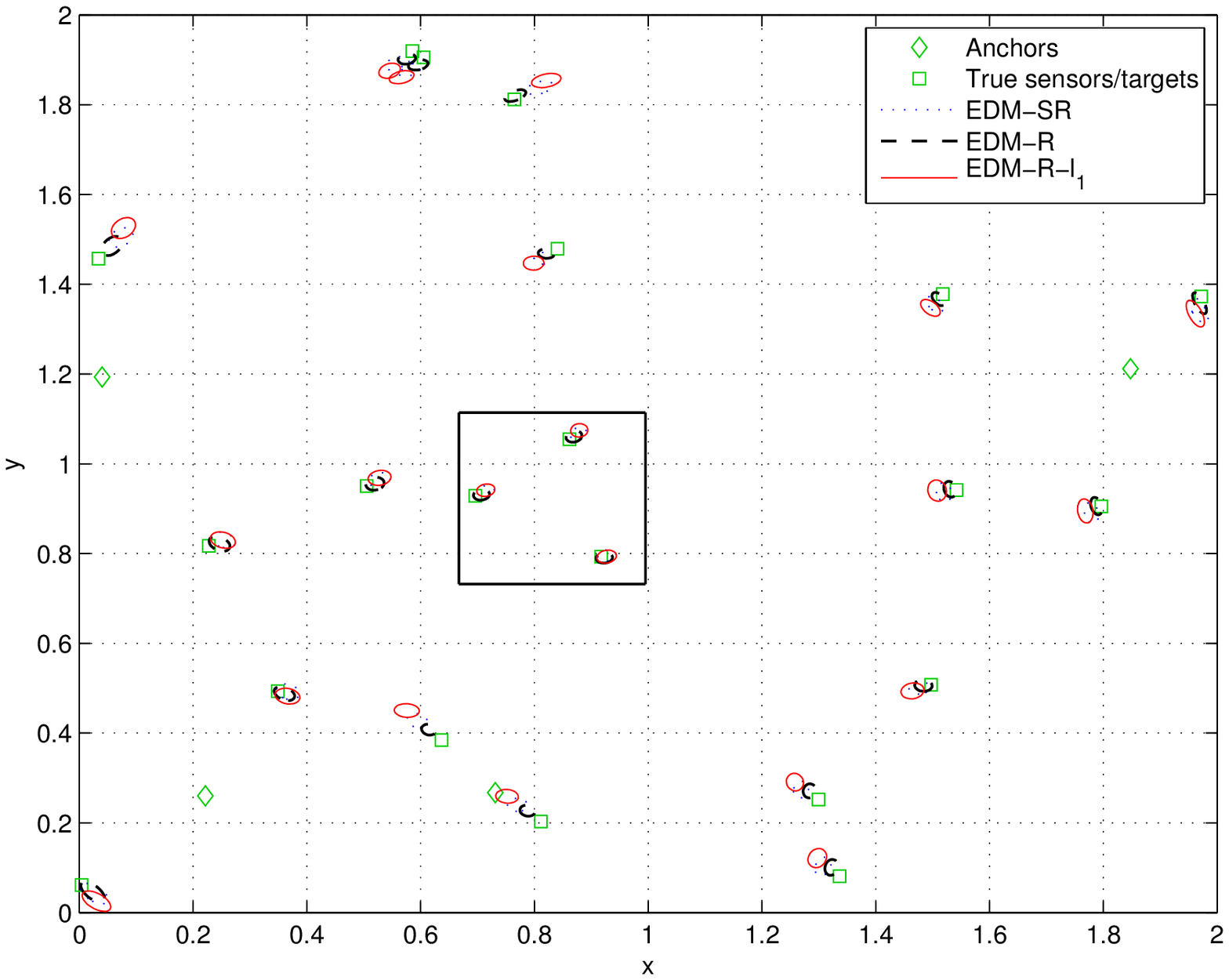}}\\
\subfloat[\label{mean_var_zoomview}Detail view]{\includegraphics[width=5.5in, height=3in]{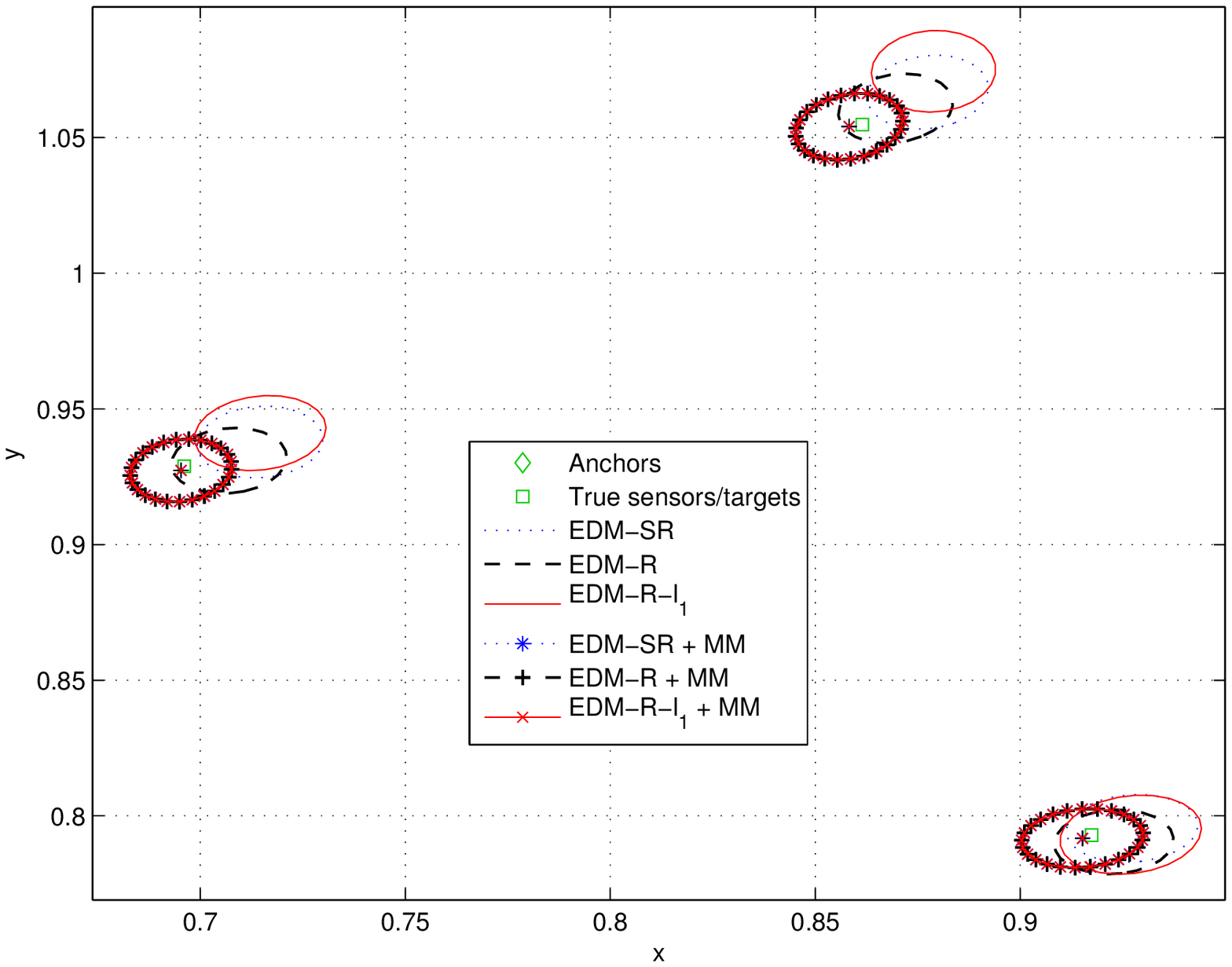}}
\end{center}
\caption{Mean and uncertainity ellipsoids of MM with different initialization methods and no outliers for $\sigma_{\text{gaussian}} = 0.025$.}
\label{mean_var}
\end{figure}

\begin{figure}%[h]
\begin{center}
\subfloat[Whole constellation]{\includegraphics[width=5.5in, height=3in]{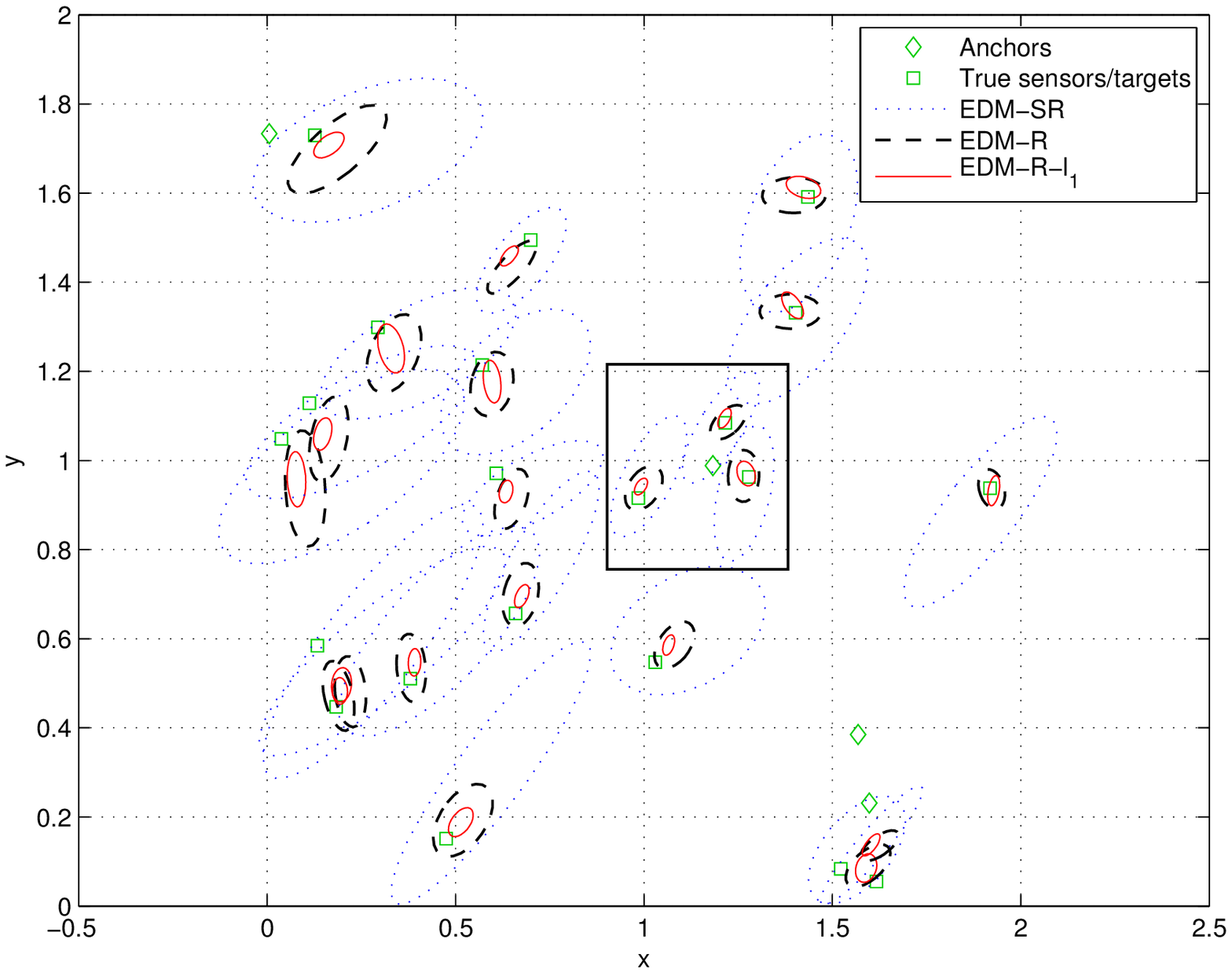}}\\
\subfloat[Detail view (including an anchor at $(1.18; \,  0.99)$)]{\includegraphics[width=5.5in, height=3in]{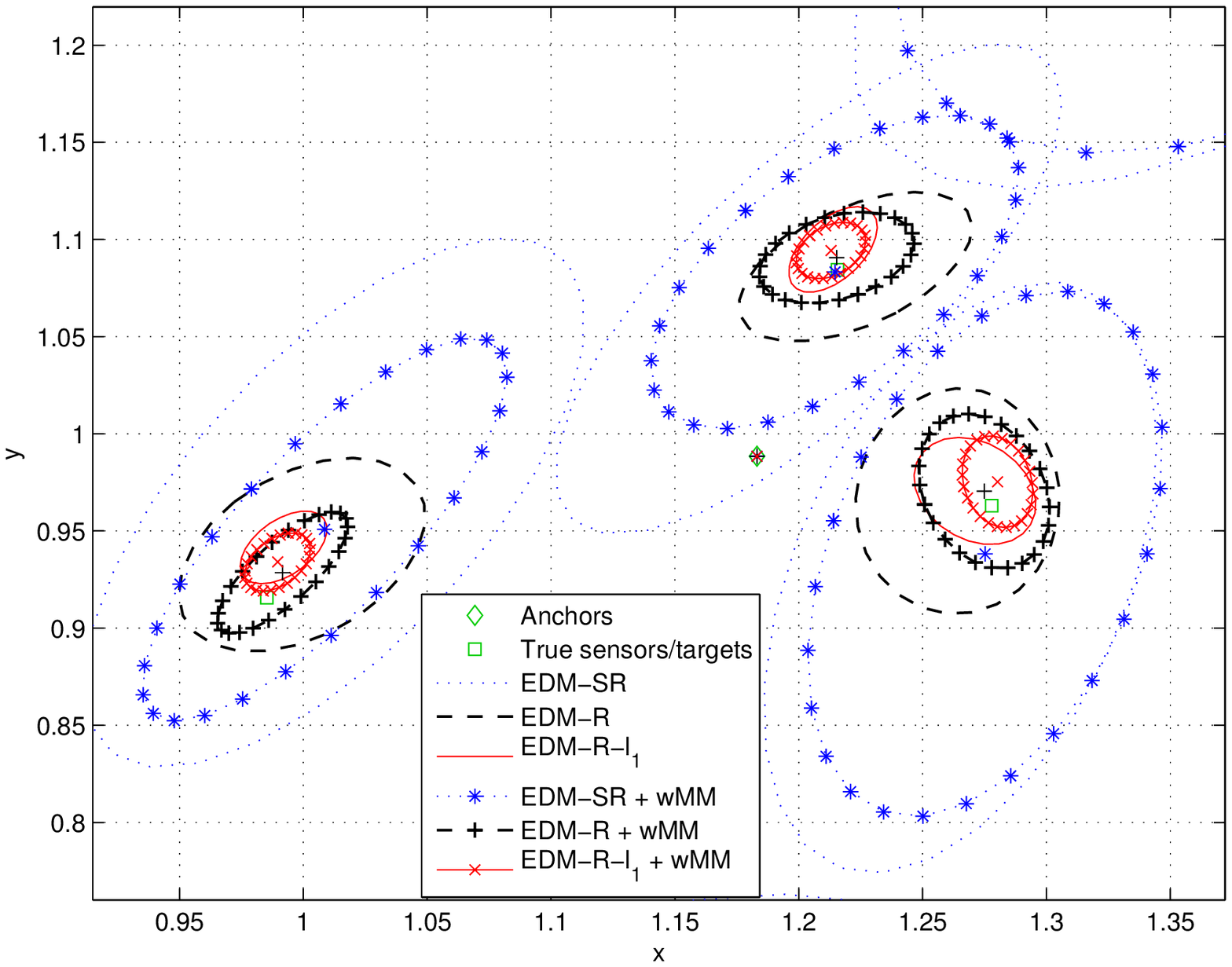}}
\end{center}
\caption{Mean and uncertainity ellipsoids of weighted MM with different initialization methods and outliers, $\sigma_{\text{outlier}} = 0.8/ \sigma_{\text{gaussian}} = 0.02$.}
\label{mean_var_outlier}
\end{figure}

\paragraph*{\textbf{Example 2 [Uncertainity Ellipsoids]}} To further examine the accuracy of MM and weighted-MM with different initialization methods, we randomly generated two networks of 10 sensors, 4 anchors and 11 target positions. 100 Monte Carlo runs were used to find the mean and ($1\sigma$) uncertainity ellipsoids of the positions estimated by the methods. The mean and uncertainity ellipsoids for $\sigma_{\text{gaussian}} = 0.025 $ and $\sigma_{\text{gaussian}} = 0.02$/$\sigma_{\text{outlier}} = 0.8$ are shown in Figs. \ref{mean_var} and \ref{mean_var_outlier}, respectively. Again, outliers are randomly added to two range measurements in Fig. \ref{mean_var_outlier}.

Without outliers (Fig. \ref{mean_var}) using EDM-SR, EDM-R, or EDM-R-$l_{1}$ as an initialization to MM makes the uncertainity ellipsoids shrink dramatically after refinement, yielding very similar means and covariances. These are only displayed in the detail view of Fig.\ \ref{mean_var_zoomview}, as they are too small to be shown in Fig.\ \ref{mean_var_wholeview}. In the presence of outliers, (Fig. \ref{mean_var_outlier}) the uncertainity ellipsoids of EDM-SR+wMM are bigger than for other methods and the means of the estimated positions are shifted. Since EDM-R-$l_{1}$ and EDM-R initializations converge to global extrema most of the time, the means of the positions estimated by weighted MM still approach the true positions and their uncertainity ellipsoids are much smaller than for EDM-SR+wMM. In the presence of outliers this example shows that EDM-R-$l_{1}$+wMM is clearly superior to the other methods.
%%%%%%%%%%%%%%%%%%%%%

\paragraph*{A Note on Practical Computational Complexity of EDM Initialization} Our experiments were conducted on a machine with an Intel Xeon 2.93 GHz Quad-Core CPU and 8 GB of RAM, using Matlab 7.1, CVX 1.2 and Yalmip 3/SeDuMi 1.1 as a general-purpose SDP solver. CPU times are similar for EDM-SR, EDM-R and EDM-R-l1, under 5 seconds for the example described above with $n = 25$ unknown positions and empirically increasing with $n^{4.5}$ for larger values of $n$ ($<100$). This gives a notion of what network sizes are currently practical for the EDM initialization methods, while keeping in mind that CPU times are known to be  unreliable surrogates for intrinsic computational complexity due to dependencies on factors such as machine hardware architecture, operating system, efficiency of numerical libraries, and solver preprocessing. No attempt was made to formulate the EDM completion problems in the most efficient way possible for the SDP solver. For MM type iterative algorithms extremely large problem sizes can be efficiently handled using contemporary numerical algorithms and computing platforms. In our experiments each iteration takes up to about 1 millisecond.
%%%%%%%%%%%%%%%%%%%%%%%%

\paragraph*{\textbf{Example 3 [Comparison of Source Localization Methods for Time-Recursive Initialization]}}  In this example SL$\ell_1$ and SLCP are compared using five anchors. We performed 100 Monte Carlo runs, where in each run the anchor and source positions were randomly generated from a uniform distribution over the square $[-10,10]\times[-10,10]$. Table \ref{tab:mse} lists the RMSE of source positions under Gaussian noise, with standard deviations $\sigma_{\text{gaussian}} = 10^{-3}$, $10^{-2}$, $10^{-1}$, and 1. SL$\ell_1$ uses $\sigma = 10^6$ for the ``projector'' $\bm{\Pi} = (\bm{\Lambda} + \sigma\mathbf{11}^T)^{-1}$, which is a very conservative value (see Appendix \ref{sec_sll1proofs}).
\begin{table}[htb]
\centering
\begin{tabular}{|c|c|c|}\hline \hline
$\sigma_{\text{gaussian}}$&SL$\ell_1$&SLCP\\ \hline
1e-3 &1.5e-3 &1.3e-3   \\
1e-2 &1.41e-2&1e-2     \\
1e-1 &0.1720 &0.1179  \\
1    &1.7922 &1.4765     \\ \hline \hline
\end{tabular}
\caption{RMSE of SL$\ell_1$ and SLCP under Gaussian noise.}
\label{tab:mse}
\end{table}

To compare the algorithms in the presence of outliers, range measurements are created according to $d_{i} = \|\cdot\|+ v_{i}$, where $v_{i}$ is a Laplacian noise term with standard deviation $\sigma_{\text{laplace}}\in [0.2,1.6]$. Results are also presented for the alternative selective Gaussian outlier generating model of Examples 1--2 with $\sigma_{\text{outlier}}\in [0.5,2]$ and $\sigma_{\text{gaussian}}=0.04$, where outliers only affect measured ranges between the second anchor and the source. Tables \ref{tab:mse_outlier:laplace} and \ref{tab:mse_outlier:selectgauss} list the RMSE of source positions for these two models.
\begin{table}[htb]
\centering
\subfloat[Laplacian noise\label{tab:mse_outlier:laplace}]{
\begin{tabular}{|c|c|c|}\hline \hline
$\sigma_{\text{laplace}}$ & SL$\ell_1$& SLCP\\ \hline
0.2 &0.2742    &0.2757\\
0.4 &0.3990    &0.4367\\
0.8 &0.8749    &0.9781\\
1.6 &1.5703    &1.8872\\ \hline \hline
\end{tabular}}
\hspace{2em}
\subfloat[Selective Gaussian noise\label{tab:mse_outlier:selectgauss}]{
\begin{tabular}{|c|c|c|}\hline \hline
$\sigma_{\text{outlier}}$ & SL$\ell_1$& SLCP\\ \hline
0.5 &0.2124    &0.2423\\
1.0 &0.3072    &0.4734\\
1.5 &0.7067    &0.7778\\
2.0 &1.0037    &1.4182\\ \hline \hline
\end{tabular}}
\caption{RMSE of SL$\ell_1$ and SLCP in the presence of outliers.}
\label{tab:mse_outlier}
\end{table}
We conclude from Tables \ref{tab:mse} and \ref{tab:mse_outlier} that the relative accuracy of SL$\ell_1$ and SLCP depends on the data generation model. For Gaussian noise the RMSE of SL$\ell_1$ is about 30\% higher than that of SLCP, whereas in the presence of outliers the situation is reversed, and SLCP exhibits an excess RMSE of 10--30\%. Interestingly, the performance gap is actually larger for selective Gaussian outliers, whose generating model does not match the assumptions of SL$\ell_1$. Similarly to Fig.\ \ref{compinit:outliers} it was found that the differences in initialization accuracy using SLCP or SL$\ell_1$ are large enough to prevent closing of the RMSE gap after weighted MM refinement due to convergence to undesirable extrema of the likelihood function.
%%%%%%%%%%%%%%%%%%%%%%%%

\paragraph*{\textbf{Example 4 [Time Recursive Updating]}} This example assesses the performance of the full time-recursive procedure, comprising SLCP or SL$\ell_1$ initialization followed by refinement. The network scenario has 16 unknown sensors, 4 anchors and 10 target locations, all randomly positioned. A new target sighting (the 11th one) becomes available and is processed incrementally, i.e., the position is estimated through SLCP or SL$\ell_1$ by fixing all the remaining ones, then all estimates are jointly refined. Results are benchmarked against refinement with full batch initialization, which makes a fresh start to the process without using any previous knowledge at every new target position to be estimated, solving different and increasingly large EDM completion problems for ML initialization.

\begin{figure}%[h]
\centering
\includegraphics[width=5.5in, height=3in]{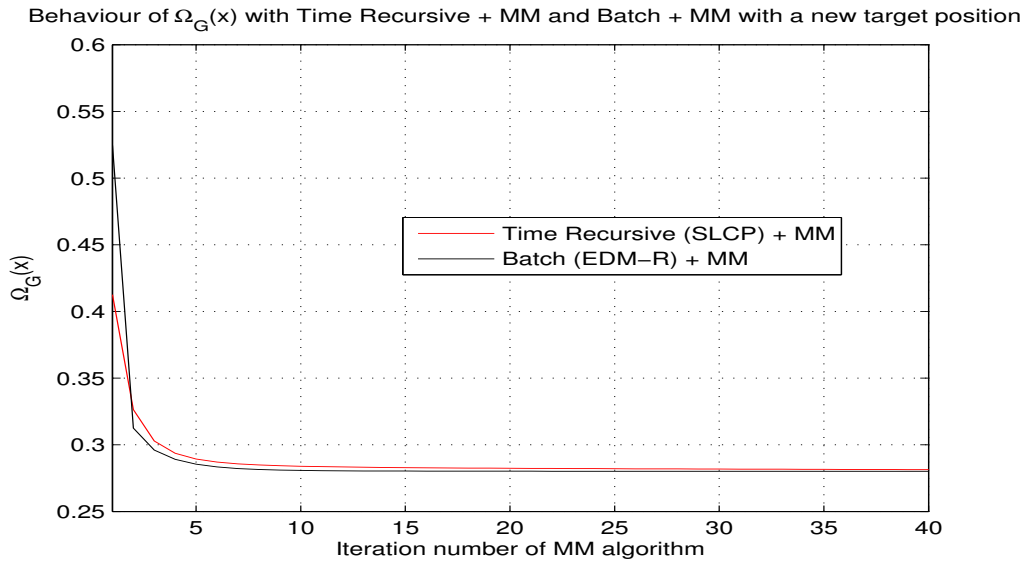}
\caption{Evolution of Gaussian cost function $\Omega_G({\bf x})$ during refinement for EDM-R+MM and SLCP+MM approaches, with $\sigma_{\text{gaussian}} = 0.04$.}
\label{dec_time_init_gaussian}
\end{figure}

This type of incremental approach was used in \cite{Ekim2009} with the SR-LS algorithm of \cite{Stoica2008} and MM refinement for Gaussian noise. SLCP is used here instead of SR-LS because, as shown in \cite{Ekim2010}, it increases the convergence speed of subsequent iterative methods and also alleviates the problem of convergence to local extrema of the ML cost function by providing better initial points than SR-LS does. Figure\ \ref{dec_time_init_gaussian} shows the evolution of the Gaussian cost function $\Omega_G({\bf x})$ during refinement after ranges to the 11th target position are sensed ($\sigma_{\text{gaussian}} = 0.04$). The time recursive (SLCP)+MM approach takes advantage of previously estimated positions to start with a lower cost than batch (EDM-R)+MM, but reaches the same final error value.       

\begin{figure}%[h]
\centering
\includegraphics[width=5.5in, height=3in]{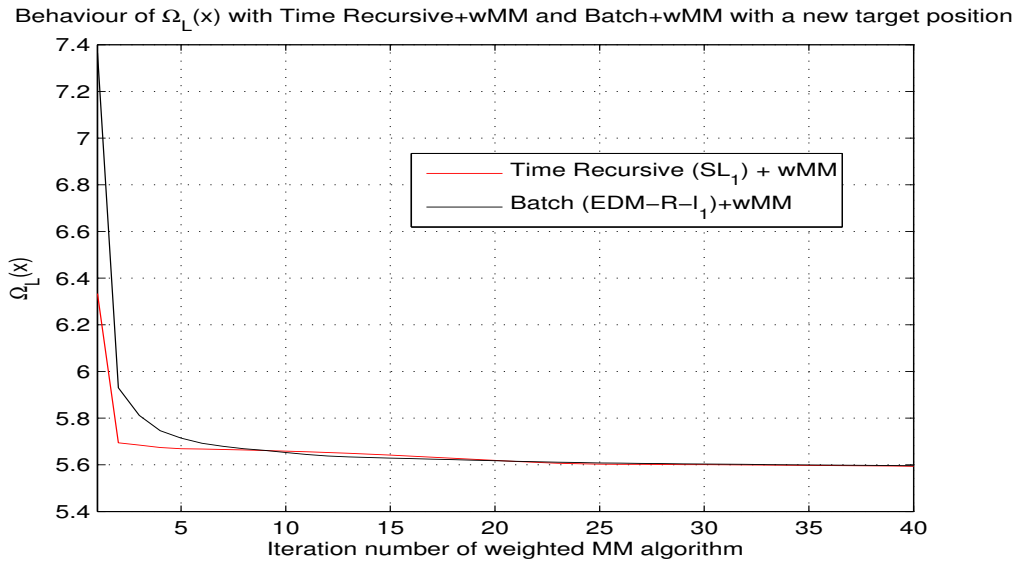}
\caption{Evolution of Laplacian cost function $\Omega_L({\bf x})$ during refinement for EDM-R-$\ell_1$+wMM and SL$\ell_1$+wMM approaches, with $\sigma_{\text{laplacian}} = 0.1$.}
\label{dec_time_init}
\end{figure}

The same network scenario is adopted in the presence of outliers. Figure\ \ref{dec_time_init} shows the evolution of cost function $\Omega_L({\bf x})$ during refinement for Laplacian outliers ($\sigma_{\text{laplacian}} = 0.1$), whose behavior is similar to the Gaussian case of Fig.\ \ref{dec_time_init_gaussian}. In both Gaussian and Laplacian settings refinement yields similar accuracy and convergence speed after batch or time-recursive initializations. Therefore, time-recursive updating is seen to retain the essential features of our EDM-based approach to SLAT, namely, very limited need for \emph{a priori} spatial information and fast convergence, at a fraction of the computational cost.

%%%%%%%%%%%%%%%%%%%%%%%
\section{Conclusion}\label{conc}
%%%%%%%%%%%%%%%%%%%%%%%%%%%%%%%%%%%%%%%%%%%%%%%%%%%%%%%%%%%%%%%%%%
In this paper, we have presented a ML-based technique to solve a SLAT problem under Gaussian and Laplacian noise. An MM method is proposed to iteratively maximize the non-convex likelihood function, for which a good initial point is required. Therefore, we have investigated two initialization schemes, namely, the batch (EDM) and time-recursive (SLCP/SL$\ell_1$) approaches that bypass the need for priors on sensor/target locations. After the first block of measurements is obtained, an EDM completion method is used for the first initialization of the sensor network topology. In our experiments this was accomplished reasonably fast (a few seconds) for scenarios with up to about 30 unknown positions. As EDM completion is not scalable, we resort to an alternative, lightweight, incremental initialization scheme as additional target range measurements become available. The SLCP or SL$\ell_1$ time recursive methods fix the already estimated positions whenever a new position is to be determined; afterwards all positions are given as initialization to the likelihood optimization methods.

Simulation results showed that our method nearly attains the Cram\'{e}r-Rao lower bound under moderate Gaussian noise. In the presence of outliers, EDM-R-$\ell_1$ or SL$\ell_1$ provide more accurate initial position estimates than other existing methods. Moreover, when used as input to iterative refinement methods they provide a good starting point that reduces the probability of convergence to undesirable extrema, yielding improved overall estimation performance. Hence, with this methodology, we obtain a processing structure that is robust to outliers and provides a scalable and accurate solution to the SLAT problem. Importantly, the algorithms based on $\ell_1$ norm optimization exhibited this robust behavior in simulation not only for Laplacian outliers, but also for an alternative outlier generation technique that did not match the underlying Laplacian modeling assumptions.

%%%%%%%%%%%%%%%%%%%%%%%%%%%%%%%%%%%%%%%%%%%%%%%%%%%%%%%%%%%%%%%%%%
\appendices
\section{Convergence of Weighted Majorization Minimization}
\label{sec_wmmconv}
To prove (local) convergence of the weighted MM iteration \eqref{itwMM} the Laplacian cost function \eqref{costfuncnew} is first majorized at time $t$ by
\begin{equation}\label{midmajfunc}
    \Gamma_{L}^{t}({\bf x}) =  \frac{1}{2}\sum_{i,j} \Bigl\{ u_{ij}^{t}(f_{ij}({\bf x})-d_{ij})^2+ \frac{1}{u_{ij}^{t}} \Bigr\} 
     + \frac{1}{2}\sum_{k,j} \Bigl\{ v_{kj}^{t}(g_{kj}({\bf x})-d_{kj})^2 + \frac{1}{v_{kj}^{t}} \Bigr\} ,
\end{equation}
where $f_{ij}$, $g_{kj}$ and $u_{ij}^{t}$, $v_{kj}^{t}$ are defined in \eqref{eq_dist} and \eqref{eq_wmmweights}. The inequality $\Omega_{L}(\mathbf{x}) \leq \Gamma_{L}^{t}(\mathbf{x})$ follows from
\begin{equation}
  \begin{split}
    \Gamma_{L}^{t}(\mathbf{x}) - \Omega_{L}(\mathbf{x}) & = \frac{1}{2}\sum_{i,j} \Bigl\{ u_{ij}^{t}(f_{ij}({\bf x})-d_{ij})^2+ \frac{1}{u_{ij}^{t}} - 2\abs{f_{ij}(\mathbf{x})-d_{ij}} \Bigr\} \\
& \quad + \frac{1}{2}\sum_{k,j} \Bigl\{ v_{kj}^{t}(g_{kj}({\bf x})-d_{kj})^2+ \frac{1}{v_{kj}^{t}} - 2\abs{g_{kj}(\mathbf{x})-d_{kj}} \Bigr\} \\
& = \frac{1}{2}\sum_{i,j} \Bigl\{ \sqrt{u_{ij}^{t}}\abs{f_{ij}({\bf x})-d_{ij}} - \frac{1}{\sqrt{u_{ij}^{t}}} \Bigr\}^{2} 
 \quad + \frac{1}{2}\sum_{k,j} \Bigl\{ \sqrt{v_{kj}^{t}}\abs{g_{kj}({\bf x})-d_{kj}} - \frac{1}{\sqrt{v_{kj}^{t}}} \Bigr\}^{2} \geq 0.
 \end{split}
\end{equation}
It is easy to check that $\Omega_{L}(\mathbf{x}^{t}) = \Gamma_{L}^{t}(\mathbf{x}^{t})$, so $\Gamma_{L}^{t}(\mathbf{x})$ has the properties of a true majorization function for the iterate $\mathbf{x}^{t}$. Now the same technique used in \eqref{expfunc} is applied to majorize \eqref{midmajfunc} by a convex quadratic function of $\mathbf{x}$, yielding
\begin{equation}\label{eq_expwfunc}
  \begin{split}  
   \Omega_{L}(\mathbf{x}) \leq \frac{1}{2}\sum_{i,j} \Bigl\{ u_{ij}^{t} \bigl( f_{ij}^{2}({\bf x})-2d_{ij}f_{ij}({\bf x}^{t}) - 2d_{ij} \langle \nabla f_{ij}({\bf x}^{t}),({\bf x}-{\bf x}^{t})\rangle +d_{ij}^{2} \bigr) + \frac{1}{u_{ij}^{t}} \Bigr\} \\
    + \frac{1}{2}\sum_{k,j} \Bigl\{ v_{ij}^{t} \bigl( g_{kj}^{2}({\bf x})-2d_{kj}g_{kj}({\bf x}^{t})-2d_{kj} \langle \nabla g_{kj}({\bf x}^{t}), ({\bf x}-{\bf x}^{t}) \rangle +d_{kj}^{2} \bigr) + \frac{1}{v_{ij}^{t}} \Bigr\}.
 \end{split}
\end{equation}
As before, equality holds for $\mathbf{x} = \mathbf{x}^{t}$, so the right-hand side of \eqref{eq_expwfunc} is still a valid majorization function. Discarding constant terms the weighted MM iteration \eqref{itwMM} results.

%%%%%%%%%%%%%%%%%%%
\section{Properties of Single-Source Localization using SL$\ell_1$}
\label{sec_sll1proofs}

\begin{IEEEproof}[Proof of Lemma \ref{thm:1}]
To streamline the notation we define $K_i = |\|y-b_i\|-d_{i}|$, and apply the KKT condition to the inner optimization problem in \eqref{sourcecostfuncl1.3} while fixing $y$. The Lagrangian function is 
\begin{equation}
L(\bm{\lambda}, \gamma) =  \sum_{i=1}^{n+l} \frac{K_i^2}{\lambda_i} + \gamma(\mathbf{1}^T\bm{\lambda}-1).
\end{equation}
The KKT conditions are
\begin{align}\label{KKT}
\frac{d\,L}{d\lambda_i} & = -\frac{K_i^2}{\lambda_i^{2}} + \gamma^* = 0, & 
\mathbf{1}^T\bm{\lambda} & = 1.
\end{align}
Using \eqref{KKT}, we find $\lambda_i^* = \frac{K_i}{\sum_{i=1}^{n+l}K_i}$ as a solution of the inner optimization problem. Plugging the optimal $\bm{\lambda}$ in the cost function of \eqref{sourcecostfuncl1.3} yields $(\sum_i K_i)^2$, thus establishing the equivalence with \eqref{sourcecostfuncl1.2}. 
\end{IEEEproof}

\begin{IEEEproof}[Approximation accuracy of $\bm{\Pi} = \lim_{\sigma \to \infty}(\bm{\Lambda} + \sigma \mathbf{1} \mathbf{1}^{T})^{-1}$]
To decide how large $\sigma$ should be, let us first define $\bm{\Pi}(\sigma) = (\bm{\Lambda} + \sigma \mathbf{1} \mathbf{1}^{T})^{-1}$. The norm of the difference to the original definition of $\bm{\Pi}$ in \eqref{sll1_projector} is given by
\begin{equation}
\begin{split}
\|\bm{\Pi} - \bm{\Pi}(\sigma)\|_{F} & = \|\bm{\Lambda}^{-1}\mathbf{1} [(\mathbf{1}^{T}\bm{\Lambda}^{-1}\mathbf{1})^{-1} - (\mathbf{1}^{T}\bm{\Lambda}^{-1}\mathbf{1} + \sigma^{-1})^{-1}] \mathbf{1}^{T}\bm{\Lambda}^{-1}\|_{F} \\
& = \frac{\mathbf{1}^{T}\bm{\Lambda}^{-2}\mathbf{1}}{(\mathbf{1}^{T}\bm{\Lambda}^{-1}\mathbf{1})(\sigma \mathbf{1}^{T}\bm{\Lambda}^{-1}\mathbf{1} + 1)}.
\end{split}
\end{equation}
Now assume the most unfavorable case with identical $\lambda_{i} = \frac{1}{n+l}$, such that
\begin{equation}
\|\bm{\Pi} - \bm{\Pi}(\sigma)\|_{F} = \frac{n+l}{\sigma(n+l)^2+1} \leq \epsilon \; \Rightarrow \; \sigma \geq \frac{1}{(n+l) \epsilon} - \frac{1}{(n+l)^2}.
\end{equation}
For $\epsilon = 10^{-4}$ and $n+l = 100$, for example, this yields $\sigma \geq 10^2-10^{-4} \approx 10^2$, which is quite low and does not raise any numerical issues in commonly available convex optimization solvers.
\end{IEEEproof}

%%%%%%%%%%%%%%%%%%%
\section{Derivation of CRLB for Gaussian Noise}
\label{sec_crlb}
The log of the joint conditional pdf for the SLAT problem is (up to an additive constant)
\begin{equation}
\label{costfunc1}
\mbox{log}f({\bf d}|{\bf x}) = -\frac{1}{2\sigma^2}\left\{\sum_{i,j}(\|x_{i}-e_{j}\|-d_{ij})^{2} + \sum_{k,j}(\|a_{k}-e_{j}\|-d_{kj})^{2}\right\},
\end{equation}
where, similarly to $\mathbf{x}$, $\mathbf{d}$ denotes the concatenation of all range measurements. Let us define matrices $\mathbf{M}_{ij}$ and $\mathbf{N}_{j}$ that extract individual positions or their differences from the vector of concatenated coordinates $\mathbf{x}$ as follows
\begin{align}
\mathbf{M}_{ij}\mathbf{x} & = x_{i} - e_{j}, & \mathbf{N}_{j}\mathbf{x} = -e_{j}.
\end{align}
Thus, (\ref{costfunc1}) is rewritten as
\begin{align}
\begin{split}
\label{ncostfunc}
\log f({\bf d}|{\bf x}) = -\frac{1}{2\sigma^2}\left\{\sum_{i,j}(\|\mathbf{M}_{ij}{\bf x}\|-d_{ij})^{2} + \sum_{k,j}(\|a_{k}+\mathbf{N}_{j}{\bf x}\|-d_{kj})^{2}\right\}.
\end{split}
\end{align}
The first derivative of (\ref{ncostfunc}) with respect to ${\bf x}$ is
%%%%%%%%%%%%%%%%
\begin{multline}
\nabla_{\bf x} \log f({\bf d}|{\bf x})  = -\frac{1}{\sigma^2}\left\{\sum_{i,j}(\|\mathbf{M}_{ij}{\bf x}\|-d_{ij})\frac{\mathbf{M}_{ij}^T\mathbf{M}_{ij}{\bf x}}{\|\mathbf{M}_{ij}{\bf x}\|}%\right.\\
%\left.
+\sum_{k,j}(\|a_{k}+\mathbf{N}_{j}{\bf x}\|-d_{kj})\frac{\mathbf{N}_{j}^T(a_{k}+\mathbf{N}_{j}{\bf x})}{\|a_{k}+\mathbf{N}_{j}{\bf x}\|}\right\}.
\end{multline}
%%%%%%%%%%%%%%%%
The second derivative of (\ref{ncostfunc}) with respect to ${\bf x}$ is
%%%%%%%%%%%%%%%
\begin{multline}
\label{2ndder}
\nabla_{\bf x}^{2} \log f({\bf d}|{\bf x})= -\frac{1}{\sigma^2}\left\{\sum_{i,j}\left\{\frac{\mathbf{M}_{ij}^T\mathbf{M}_{ij}{\bf x}{\bf x}^T\mathbf{M}_{ij}^T \mathbf{M}_{ij}}{\|\mathbf{M}_{ij}{\bf x}\|^2} %\right.\right. \\
%\left.
+\frac{\|\mathbf{M}_{ij}{\bf x}\|-d_{ij}}{\|\mathbf{M}_{ij}{\bf x}\|}
\left(\mathbf{M}_{ij}^T\mathbf{M}_{ij}-\frac{\mathbf{M}_{ij}^T\mathbf{M}_{ij}{\bf x}{\bf x}^T\mathbf{M}_{ij}^T\mathbf{M}_{ij}}{\|\mathbf{M}_{ij}{\bf x}\|^2}\right)\right\}\right. \\
\left.+\sum_{k,j}\left\{\frac{\mathbf{N}_{j}^T(a_{k}+\mathbf{N}_{j}{\bf x})(a_{k}+\mathbf{N}_{j}{\bf x})^T\mathbf{N}_{j}^T}{\|a_k+\mathbf{N}_{j}{\bf x}\|^2}%\right.\\
%\left.\left.
+\frac{\|a_k+\mathbf{N}_{j}{\bf x}\|-d_{kj}}{\|a_k+\mathbf{N}_{j}{\bf x}\|}\left( \mathbf{N}_{j}^T\mathbf{N}_{j}
-\frac{\mathbf{N}_{j}^T(a_k+\mathbf{N}_{j}{\bf x})(a_k+\mathbf{N}_{j}{\bf x})^T\mathbf{N}_{j}^T}{\|a_k+\mathbf{N}_{j}{\bf x}\|^2} \right)
\right\}\right\}
\end{multline}
%%%%%%%%%%%%%%%%%
The Fisher information matrix, $F_{\bf x}$, is obtained by taking the negative expected value of the (\ref{2ndder}) with respect to ranges as \cite{Kay1993}
%%%%%%%%%%%%%%%%%%
\begin{multline}
F_{\bf x} = -\mbox{E}_{{\bf d}}\{\nabla_{\bf x}^{2} \log f({\bf d}|{\bf x})\}=\frac{1}{\sigma^2}\left\{\sum_{i,j}\frac{\mathbf{M}_{ij}^T\mathbf{M}_{ij}{\bf x}{\bf x}^T\mathbf{M}_{ij}^T \mathbf{M}_{ij}}{\|\mathbf{M}_{ij}{\bf x}\|^2}%\right.\\
%\left.
+\sum_{k,j}\frac{\mathbf{N}_{j}^T(a_{k}+\mathbf{N}_{j}{\bf x})(a_{k}+\mathbf{N}_{j}{\bf x})^T\mathbf{N}_{j}^T}{\|a_k+\mathbf{N}_{j}{\bf x}\|^2}\right\}.
\end{multline}
%%%%%%%%%%%%%%%%%
The CRLB matrix in (\ref{crbound}) is taken as the inverse of $F_{\bf x}$.

\bibliographystyle{IEEEtran}
\bibliography{strings1}

%\begin{IEEEbiographynophoto}{P\i nar~O\u{g}uz-Ekim}
%\begin{IEEEbiographynophoto}{Jo\~{a}o~Gomes}
%\begin{IEEEbiographynophoto}{Jo\~{a}o~Xavier}
%\begin{IEEEbiographynophoto}{Paulo~Oliveira}

\end{document}